\documentclass{article}
\usepackage[a4paper, total={5.5in, 8in}]{geometry}
\usepackage{dutchcal}

\newcommand{\mr}[1]{\mathring{#1}}
\usepackage{graphics}
  \usepackage{amssymb}
    \usepackage{amsthm}
  \usepackage{amsmath}
  \usepackage{pgf}
\usepackage{tikz}
\usetikzlibrary{arrows,automata}
\usepackage{epstopdf}
  
\def\p{\partial}

\newtheorem{theorem}{Theorem}[section]
\newtheorem{lemma}[theorem]{Lemma}

\newtheorem{proposition}[theorem]{Proposition}
\theoremstyle{definition}

\newtheorem{example}[theorem]{Example}
\newtheorem{remark}[theorem]{Remark}

\newcommand {\A}{\mathbb{A}}
\newcommand {\bay}{\begin{array}}
\newcommand {\eay}{\end{array}}
\newcommand {\bdm}{\begin{displaymath}}
\newcommand {\edm}{\end{displaymath}}
\newcommand{\ups}{\upsilon}
\newcommand{\w}{\varpi}

  \makeatletter
  \let\c@equation\undefined
  \let\c@section\undefined
  \let\c@subsection\undefined
  \let\c@zad\undefined
  \makeatother
  \newcounter{section}

  \newcounter{equation}[section]
  
  \newcounter{subsection}[section]

\newcommand{\beq}{\begin{equation}}
\newcommand{\eeq}{\end{equation}}
\newcommand{\R}{\mathbb R}
\newcommand{\N}{\mathbb N}

\newcommand{\bd}{\begin{displaymath}}
\newcommand{\ed}{\end{displaymath}}
\newcommand{\be}{\begin{equation}}
\newcommand{\ee}{\end{equation}}
\newcommand{\bdow}{\begin{proof}}
\newcommand{\edow}{\end{proof}}

\newcommand{\papap}{\papa \end{proof}}

\newfont{\smoldita}{cmmib8}
\newfont{\boldita}{cmmib10}
\newfont{\bboldita}{cmmib10}

\newcommand {\sem}[1]{\mbox{$({#1}(t))_{t \geq 0}$}}

\newcommand{\la}{\lambda}
\newcommand{\lam}{\lambda}
\newcommand{\mb}[1]{\boldsymbol{#1}}
\newcommand{\mbb}[1]{\mathbb{#1}}
\newcommand{\mc}[1]{\mathcal{#1}}

\newcommand{\msf}[1]{\mathsf{#1}}

\newcommand{\bv}{\mathbf{v}}
\newcommand{\Q}{\mathcal{Q}}

\begin{document}
\begin{center}{\Large Telegraph systems on networks and port-Hamiltonians. III. Explicit representation and long-term behaviour}\end{center}

\begin{center}{J. Banasiak\footnote{The research has been partially supported by the National
Science Centre of Poland Grant 2017/25/B/ST1/00051 and the National Research Foundation of South Africa Grant 82770} \\ \small{Department of Mathematics and Applied Mathematics, University of Pretoria}\\ \small{Institute of Mathematics,  \L\'{o}d\'{z} University of Technology}\\ \small{International Scientific Laboratory of
Applied Semigroup Research, South Ural
State University}\\ \small{e-mail: jacek.banasiak@up.ac.za}\\\& \\A. B\l och\footnote{The research was completed while the author was a Doctoral Candidate in the Interdisciplinary Doctoral School at \L\'{o}d\'{z} University of Technology, Poland.}\\
 \small{Institute of Mathematics,  \L\'{o}d\'{z} University of Technology} \\\small{e-mail: adam.bloch@dokt.p.lodz.pl}}\end{center}
\begin{abstract}
In this paper we present an explicit formula for the semigroup governing the solution to hyperbolic systems on a metric graph, satisfying general linear Kirchhoff's type boundary conditions. Further, we use this representation to establish the long term behaviour of the solutions. The crucial role is played by the spectral decomposition of the boundary matrix.
\end{abstract}
 \section{Introduction}
In this paper we consider systems of linear hyperbolic equations on a bounded interval, say, $[0,1]$, known also as port-Hamiltonians, \cite{JaMoZw}, coupled by boundary conditions relating the incoming and outgoing values of the solution at the endpoints $x=0$ and $x=1$. In particular, we study systems of the form
 \begin{subequations}\label{syst1}
 \begin{equation}
 \p_t\binom{\mb{\ups}}{\mb{\w}} = \left(\begin{array}{cc}-\mc C_+&0\\0&\mc C_-\end{array}\right)\partial_x\binom{\mb{\ups}}{\mb{\w}}+\mc K\binom{\mb{\ups}}{\mb{\w}},\quad 0<x<1, t>0,\label{syst1eq}
 \end{equation}
  \begin{equation}
 \mb{\ups}(x,0)= \mr{\mb\ups}(x),\; \mb{\w}(x,0)= \mr{\mb\w}(x),\quad 0<x<1,\label{syst1ic}
 \end{equation}
 \begin{equation}
 \mb \Xi_{out} (\mb{\ups}(0,t),\mb{\w}(1,t))^T + \mb\Xi_{in}(\mb{\ups}(1,t),\mb{\w}(0,t))^T=0,\quad t>0,\label{syst1bc}
 \end{equation}
 \end{subequations}
  where $\mb\ups$ and $\mb\w$ represent the densities of the flow from $0$ to $1$ and from $1$ to $0$, respectively, $\mc C_+$ and $\mc C_-$ are $m_+\times m_+$ and  $m_-\times m_-$ diagonal matrices with positive entries, $\mc K$ is a $2m\times 2m$ matrix, where $2m=m_++m_-$, $\mb \Xi_{out}, \mb\Xi_{in}$ are $2m\times 2m$ matrices relating the outgoing  $\mb{\ups}(0),\mb{\w}(1)$ and incoming $\mb{\ups}(1),\mb{\w}(0)$ flows at the boundary points.

  An important class of such problems arises from dynamical systems on metric graphs.  Let $\Gamma$ be a graph with $r$ vertices $\{\bv_j\}_{1\leq j\leq r}=:\Upsilon$ and $m$ edges $\{\mb e_j\}_{1\leq j\leq m}$ (identified with $[0,1]$ through a suitable parametrization).  The dynamics on each edge $\mb e_j$ is described by
 \begin{equation}
 \p_t \mb p^j+ \mc M^j\p_x\mb p^j+\mc N^j\mb p^j  =0, \quad t>0, 0<x<1, 1\leq j\leq m,
 \label{sys1}
 \end{equation}
 where $\mb p^j = (p^j_1, p^j_2)^T$, $\mc M^j = (M^j_{lk})_{1\leq k,l\leq 2}, \mc N^j = (N^j_{lk})_{1\leq k,l\leq 2}$ are real matrix functions defined on $[0,1]$. We assume that $\mc M^j$ are continuously differentiable and $\mc N^j$ are bounded on $[0,1]$. The central assumption is that $\mc M^j(x)$ is a strictly hyperbolic matrix for each $x\in [0,1]$ and $1\leq j\leq m$. System \eqref{sys1} is complemented with initial conditions and  suitable transmission conditions coupling the values of $\mb p^j$ at the vertices which the edges $\mb e_j$ are incident to. Then, \eqref{syst1} can be obtained from \eqref{sys1} by diagonalization so that (suitably re-indexed) $\mb \ups$ and $\mb \w$ are the Riemann invariants of $\mb p = (\mb p^j)_{1\leq j\leq m}$, see \cite[Section 1.1]{BaCor}.

 Such problems have been a subject of extensive research in the recent years. Let us mention here e.g.  \cite{AM2, DKNR, BN, BFN3, KMN, Kuch, DM} for the dynamics on graphs point of view, and \cite{BaCor, Zwart2010, JaZwbook, JaMoZw} for the 1-D hyperbolic systems point of view. Many of these papers are focused on well-posedness of the problem under various boundary conditions ensuring asymptotical stability of the resulting  semigroup. There is a fairly comprehensive theory of the long term behaviour of the solutions to transport problems on networks (which is a special case of \eqref{syst1}), but there seems to be no its counterpart for general 1-D hyperbolic systems. In this paper we try to fill this gap.

 There are two main approaches for studying the long term asymptotics of semigroups. The first and most powerful is the spectral theory for semigroups and generators, used, for instance, in \cite{KS2005} for the transport semigroup with Kirchhoff's-type boundary conditions. However, it requires many advanced tools from functional and complex analysis. The second, potentially easier, approach is to find an explicit representation of the semigroup. Though it is not always possible, if we succeed, then we can employ direct calculations or spectral theory for simpler objects such as matrices.

 An explicit formula for the transport semigroup with unit velocities appeared first in \cite{Dorn} and it was used in  \cite{BN}  for investigating long term behaviour of the semigroup solving a transport problem on a network, with more general boundary conditions than those in \cite{KS2005}. We mention that these papers focused on finding conditions ensuring the periodicity of the limit semigroup and it turned out that both methods yielded the same results, \cite{MatSik}.

 In this paper we extend the approach of \cite{BN} to general hyperbolic systems.
 We would like to emphasize that, while the considerations in \cite{KS2005,BN} are carried out in the $L^1$ setting, our theory works equally well in any $L^p$ space with $p\in[1,\infty)$ and does not depend on the value of the exponent $p$. Since the problems in \cite{KS2005,BN} fit into our framework as well,  the main result of this paper, that is, Theorem \ref{maintheorem}, is a generalization of the existing results.

  The considerations in this paper are mostly done for the principal part of \eqref{syst1}, that is, the system with $\mc K=0$. We recognize that this is a serious restriction but, as $\mc K$ induces a bounded perturbation, the structure of the full semigroup associated with \eqref{syst1} is well understood once we know the semigroup generated by its principal part by, say, the Phillips--Dyson expansion, see e.g. \cite[Theorem III.1.10]{EN}. The construction of an explicit formula for the solution to the full problem \eqref{syst1} is possible in some special cases such as unidirectional transport, see \cite[Theorem 2.9]{BaPu2019} or, in a similar way, if $\mc K$ is block diagonal with blocks corresponding to the same directions of transport. Otherwise, the approach presented here leads to problems with shifted argument, see e.g. \cite[Example 1]{Puch}.

 The paper is organized as follows. In Section \ref{earlierresults} we briefly recall the notation and results concerning the well-posedness of \eqref{syst1} from \cite{JBAB1}. Section \ref{expforsec} is focused on deriving the explicit formula.  First, we convert the problem \eqref{syst1} (with $\mc K=0$) to an equivalent one, but with unit velocities. Next, we construct the formula for the unperturbed problem  and show the relation of our semigroup to the transport semigroup.  Section \ref{asymptoticbehaviour} deals with the asymptotic behaviour of the semigroup and in the last section we present some examples illustrating our theory.

 \section{Notation, definitions and earlier results}\label{earlierresults}

   We consider a network represented by a finite, connected and simple (without loops and multiple edges) metric graph $\Gamma$ with $r$ vertices $\{\bv_j\}_{1\leq j\leq r}=:\Upsilon$ and $m$ edges $\{\mb e_j\}_{1\leq j\leq m}$. Let $E_{\bv}$ be the set of edges incident to $\bv,$ $J_{\bv} :=\{j;\; \mb e_j \in E_{\bv}\} $  and  $|E_{\bv}| =|J_{\bv}|$ be the valency of $\bv$. Each edge is identified with the unit interval through a sufficiently smooth invertible function $l_j: \mb e_j \mapsto [0,1]$. In particular, we call $\bv$ with $l_j(\bv) =0$ the tail of $\mb e_j$ and the head if $l_j(\bv)=1$.  On each edge $\mb e_j$ we consider the system \eqref{sys1}. Let $\la^j_-<\la^j_+$ be the eigenvalues of $\mc M^j, 1\leq j\leq m$ (the strict inequality is justified by the strict hyperbolicity of $\mc M^j$).  The eigenvalues can be of the same sign as well as of different signs. In the latter case, we set $\la^j_-<0<\la^j_+$. By $ f^j_\pm = (f^j_{\pm,1}, f^j_{\pm,2})^T$ we denote the eigenvectors corresponding to $\la^j_\pm$, respectively, and by $$
 \mc F^j = \left(\begin{array}{cc} f^j_{+,1}&f^j_{-,1}\\
   f^j_{+,2}&f^j_{-,2} \end{array}\right)
   $$
   the diagonalizing matrix on each edge. The Riemann invariants $\mb u^j = (u^j_1,u^j_2)^T, 1\leq j\leq m,$ are defined by
 \begin{equation}
 \mb u^j = (\mc F^{j})^{-1} \mb p^j\quad\text{and}\quad
 \mb p^j = \binom{f^j_{+,1} u^j_1 + f^j_{-,1}u^j_2}{f^j_{+,2}u^j_1 + f^j_{-,2}u^j_2}.\label{pguw}
 \end{equation}
Then, we diagonalize \eqref{sys1} as
  \begin{equation}
 \p_t\mb u^j = \left(\begin{array}{cc}-\la^j_+&0\\0&-\la^j_- \end{array}\right)\p_x\mb{u}^j+\overline{\mc N^j}\mb u^j\label{sysdiag}
 \end{equation}
 for each $1\leq j\leq m$. Our assumptions ensure that $\mb u^j\mapsto\overline{\mc N^j}\mb u^j$ induces a bounded perturbation in any $(L^p(0,1))^{2m}$ space, hence it is irrelevant for the generation of a semigroup. As noted in Introduction, further considerations are carried out with $\overline{\mc N^j}=0$.

\subsection{Boundary conditions --- from networks to port-Hamiltonians}

The most general linear local boundary conditions for \eqref{sys1} at any vertex $\bv\in \Upsilon$ can be written as
\begin{equation}
  \mb \Phi_{\bv} \mb p(\bv) = 0,
  \label{bc2s}
  \end{equation}
where $\mb p(\bv) = ((p_1^j(\bv),p_2^j(\bv))_{j\in J_{\bv}})^T$ and $\mb \Phi_{\bv}$ is a real $l\times|J_{\bv}|$ matrix, where $l$ is the number of equations relating the boundary values of $\mb p^j$s at $\bv$. Such a formulations is, however, not convenient as it does not provide a clear set of conditions on $\mb\Phi_{\bv}$ and, in particular, on $l$ that ensure the well-posedness of \eqref{sys1}, \eqref{bc2s}. To proceed, we employ the paradigm introduced in \cite[Section 1.1.5.1]{BaCor} requiring that at any vertex the outgoing data should be determined by the incoming ones. Since for \eqref{sys1} it is, in general, impossible to tell which data are outgoing and which are incoming, we re-write \eqref{bc2s} in terms of the Riemann invariants, defined by \eqref{pguw}, as
\begin{equation}
\mb \Psi_{\bv} \mb u(\bv) : = \mb \Phi_{\bv} \mc F(\bv)\mb u(\bv)  = 0,
  \label{bc2uw}
  \end{equation}
where $\mc F(\bv)=diag\{\mc F^j(\bv)\}_{j\in J_{\bv}}.$ We observe that \eqref{bc2uw} makes sense only for $\bv$ that is not a \textit{sink}, that is, a vertex with no outgoing data, see \cite[Definition 2.3]{JBAB2}. The boundary condition \eqref{bc2uw} is called  \textit{generalized Kirchhoff's condition} at $\bv$.

For that, we notice that by the continuity and strict hyperbolicity of $\mc M^j, 1\leq j\leq m$, the eigenvalues $\la_+^j, \la_-^j$ are never zero and hence each equation in \eqref{sysdiag} describes a flow in a fixed direction. Hence, we only need to distinguish functions describing the flow from $0$ to $1$ and from $1$ to $0$, with no reference to the network structure.  Accordingly, we group the Riemann invariants $\mb u$  into parts corresponding to positive and negative eigenvalues and rename them as
 \begin{equation}
\begin{split}
\mb {\upsilon}  &:=\left((u^j_1)_{j\in J_1\cup J_2},(u^j_2)_{j\in J_2}\right) = (\ups_j)_{j\in J^+},\\
\mb {\w}& := \left((u^j_1)_{j\in J_0}, (u^j_2)_{j\in J_1\cup J_0}\right) = (\w_j)_{j\in J^-},
\end{split}
\label{renum}
\end{equation}
where $J^+$  and $J^-$  are the sets of indices $j$ with at least 1 positive eigenvalue, and at least 1 negative eigenvalue of $\mc M^j$, respectively. Since in  $J^+$ (respectively $J^-$) the indices from $J_2$ (respectively $J_0$) appear twice, we renumber them in some arbitrary (but consistent) way to avoid confusion. This results in just re-labelling of the components of \eqref{syst1} without changing its structure.\\
This procedure converts the problem \eqref{sys1} on $\Gamma$ into a first order transport problem on a multi digraph $\mb\Gamma$ with the same vertices $\Upsilon$ and where each edge in $\Gamma$ was split into two edges in $\mb\Gamma$ paramterized by $x\in[0,1]$. Then, after combining the matrices $\mb \Psi_\bv$ over all vertices $\bv$ that are not sinks to a global matrix $\mb\Psi,$ splitting the latter into the outgoing and incoming parts and re-indexing, see \cite{JBAB1}, the boundary conditions \eqref{bc2uw} can be written as
$$
\mb\Xi(\mb\upsilon(0),\mb\upsilon(1),\mb\w(0),\mb\w(1))^T=\mb\Xi_{out}(\mb\upsilon(0),\mb\w(1))^T+\mb\Xi_{in}(\mb\upsilon(1),\mb\w(0))^T=0.
$$
This formulation does not depend on the fact that $\mb\Xi$ has a special form coming from  Kirchhoff's boundary conditions but it can be an arbitrary $2m\times 4m$ matrix. Hence, we arrive at the problem \eqref{syst1}, where
$$
\mc C_+=diag(c_j(x))_{j\in J^+},\quad \mc C_-=diag(c_j(x))_{j\in J^-}
$$
and the functions $c_j$ equal the absolute values of the corresponding eigenvalues.
\subsection{Well-posedness of \eqref{syst1}}
Without loss of generality, see \cite[Theorem 3.3 (2)]{Zwart2010}, we assume that $\mb\Xi_{out}$ is invertible and introduce the matrix $\mc B:=\mb\Xi_{out}^{-1}\mb\Xi_{in}$. Finally, we denote $\mb X_p:=(L^p(0,1))^{2m}$ and define the operator $(A_{p,\mc B}, D(A_{p,\mc B}))$ in $\mb X_p$ as $\left.\msf A\right|_{D(A_{p,\mc B})}$, where
\begin{subequations}\label{generatordefinition}
\begin{equation}\label{gendef1}
\msf A:=diag(-\mc C_+,\mc C_-)diag(\underbrace{\partial_x,\ldots,\partial_x}_{2m\;\text{times}}),
\end{equation}
\begin{equation}\label{gendef2}
D(A_{p,\mc B}):=\left\{\binom{\mb\upsilon}{\mb\w}\in(W_1^p(0,1))^{2m}:\binom{\mb\upsilon(0)}{\mb\w(1)}=\mc B\binom{\mb\upsilon(1)}{\mb\w(0)}\right\}.
\end{equation}
\end{subequations}
Then, combining \cite[Thms. 4.1 \& 4.2]{JBAB1} we have the following
\begin{theorem}\label{genthm}
Let $\mc B$ be an arbitrary matrix. The operator $(A_{p,\mc B}, D(A_{p,\mc B}))$ generates a $C_0$-semigroup $(G_p(t))_{t\geq 0}$ on $\mb X_p$ for any $1\leq p<\infty$. Moreover, for $p>1,$ the semigroup $(G_p(t))_{t\geq 0}$ is the restriction of the semigroup $(G_1(t))_{t\geq 0}$ to the space $\mb X_p$.
\end{theorem}
The case $p=1$ was proved in \cite{Zwart2010} (see also \cite{EngKra}, \cite{JaMoZw}), however, the proof there is based on control theory methods from \cite{Staff}. In \cite{JBAB1} we provided an alternative, purely semigroup-theoretic proof based on the result from \cite{BFN3}. Then, the well-posedness for $p>1$ follows from direct estimates of the $L^1$ solutions with $L^p$ data.

\begin{remark}\label{uwagidoteorii}
Since the proof of the well-posedness uses only the form \eqref{syst1} that does not take into account any particular feature of the systems on the graph, we can make the following points.
\begin{enumerate}
\item Restricting our attention in \cite{JBAB1, JBAB2} to $2\times 2$  systems is purely for notational convenience -- the theory remains valid for systems of arbitrary (but finite) dimension on each edge.
\item The hyperbolicity of the matrices $\mc M^j(x)$ is not necessary as long as each of them is diagonalizable with nonvanishing and differentiable in $x$ eigenvalues -- e.g. when the eigenvalues are constant and semisimple.
\item The interval $[0,1]$ can be replaced by  intervals $[0,l_j],$ with $l_j$ varying between the edges.
\end{enumerate}
\end{remark}

\section{An explicit formula}\label{expforsec}
We emphasise that our theory works in each $L^p$ space and does not depend on $p$. Thus, in the remaining part of the paper, we do not make any specific assumption on $p\in[1,\infty)$.

\subsection{Conversion to unit velocities}\label{convvelsection}
Our starting point is \eqref{syst1} with $\mc K=0$ and the boundary condition \eqref{syst1bc} solved with respect to the outgoing boundary values, that is, in the form that appears in \eqref{gendef2}.

In this subsection we show that under the assumption that the traverse times along each edge are natural multiples of a one reference time,  problem \eqref{syst1} can be reduced to an equivalent problem with $c_j\equiv 1$ for all $j\in\{1,\ldots,2m\}$. The idea was first introduced in \cite{KS2005} for the transport problem with constant velocities. In \cite{BN}, the authors used this assumption to convert that problem to a transport problem with unit velocities. In \cite{BanLAA, BaPu2019}, the authors described a similar conversion for $x$-dependent velocities.  Here we present a detailed construction for the considered problem.

For a given $j\in\{1,\ldots,2m\},$ we define a function $L_j:[0,1]\to[0,L_j(1)]$ by
\bd
L_j(x)=\int\limits_0^x\frac{1}{c_j(s)}\,ds.
\ed
Then, we adopt the following assumption:
\be\label{ratdep}
\exists_{c\in\R}\;\forall_{j\in\{1,\ldots,2m\}}\quad l_j:=cL_j(1)\in\mathbb{N}.
\ee
\begin{remark}
We observe that $L_j(1)$ is the time needed to traverse the edge $\mb e_j$ with the velocity $c_j$ from the tail at $x=0$ to the head at $x=1$ or in the reversed direction (depending on whether $j\in J^+$ or $j\in J^-$), see \cite{JBAB1}. Hence, the assumption \eqref{ratdep} states that  all traverse times are natural multiples of a single reference time.
\end{remark}
By rescaling time as $\tau=ct$ and introducing a new spatial variable $y=cL_j(x)$ for each $j$, we convert the differential equations in  \eqref{syst1} to
\be
\begin{split}
\p_{\tau}\ups_j(y,\tau)=-\p_y\ups_j(y,\tau),\quad \tau>0, y\in(0,l_j),j\in J^+,\\
\p_{\tau}\w_j(y,\tau)=\p_y\w_j(y,\tau),\quad \tau>0,  y\in(0,l_j), j\in J^-,
\end{split}
\ee
see \cite[Section 2.2.5]{BaPu2019}.
In the above problem, the velocities equal 1 at the cost of different lengths of the intervals. However, due to the assumption \eqref{ratdep}, we can divide each interval $[0,l_j]$ into $l_j$ intervals of unit length, which are then translated to the interval $[0,1]$ and become the new edges. The dividing points will become the new endpoints. Then, each function is identified with an $l_j$-tuple of functions defined on the new edges, where, to preserve the structure of the original problem, we require the continuity of the functions in each $l_j$-tuple across  the new endpoints. Following the preceding discussion, we introduce the notation
\bd
\mb\nu=(\nu_j)_{j\in J^+}=((\nu_{j,i})_{i=1,\ldots,l_j})_{j\in J^+}, \mb\omega=(\omega_j)_{j\in J^-}=((\omega_{j,i})_{i=1,\ldots,l_j})_{j\in J^-}
\ed
and $\ell:=\sum_{j=1}^{2m}l_j$, which is the dimension of the new system. Since each $L_j$ is  strictly increasing, the inverse $L_j^{-1}$ exists and we can define an operator $\mb\Q:\mb X_p\to(L^p(0,1))^{\ell}$ by
\bd
\left(\mb\Q\binom{\mb\ups}{\mb\w}\right)_{j,i}(y)=\nu_{j,i}(y):=\upsilon_j\left(L_j^{-1}\left(\frac{y+i-1}{c}\right)\right), \quad j\in J^+,
\ed
\bd
\left(\mb\Q\binom{\mb\ups}{\mb\w}\right)_{j,i}(y)=\omega_{j,i}(y):=\varpi_j\left(L_j^{-1}\left(\frac{l_j+y-i}{c}\right)\right), \quad j\in J^-,
\ed
where $y\in[0,1]$, $i=1,\ldots,l_j$. The map $\mb\Q$ provides a correspondence between the original variables $(\mb\ups,\mb\w)^T$ and the new variables $(\mb\nu,\mb\omega)^T$. Precisely speaking, the function $\nu_{j,i}$ represents the values of the function $\ups_j$ on the $i$-th subinterval of the interval $[0,1]$, while the function $\omega_{j,i}$ represents the values of the function $\w_j$ on the $(l_j-i+1)$-th subinterval of $[0,1]$.
\begin{remark}
We could keep the order of identification between $\mb\w$ and $\mb\omega$ the same as for $\mb\ups$ and $\mb\nu$, that is, from $0$ to $1$. However, for $j\in J^-$, the flow described by the function $\w_j$ occurs from $1$ to $0$, thus it seems reasonable to reverse the order of the identification.
\end{remark}
The operator $\mb\Q$ transforms $\eqref{syst1}$ to the following, equivalent, problem:
\be\label{ibvpunitvel}
\begin{split}
\p_{\tau}\mb\nu(y,\tau)=-\p_y\mb\nu(y,\tau),\quad \tau>0, y\in(0,1),\\
\p_{\tau}\mb\omega(y,\tau)=\p_y\mb\omega(y,\tau),\quad \tau>0, y\in(0,1),\\
\mb\nu(y,0)=\mr{\mb\nu}(y), \mb\omega(y,0)=\mr{\mb\omega}(y),\quad y\in(0,1),\\
\binom{\mb\nu(0,\tau)}{\mb\omega(1,\tau)}=\widetilde{\mc B}\binom{\mb\nu(1,\tau)}{\mb\omega(0,\tau)},\quad \tau>0,
\end{split}
\ee
where the $\ell\times\ell$ matrix $\widetilde{\mc B}$ describes the new boundary condition for the converted problem. Precisely speaking,
\bd
\nu_{j,1}(0)=\sum\limits_{k\in J^+}b_{jk}\nu_{k,l_k}(1)+\sum\limits_{k\in J^-}b_{jk}\omega_{k,l_k}(0),\quad j\in J^+,
\ed
\bd
\omega_{j,1}(1)=\sum\limits_{k\in J^+}b_{jk}\nu_{k,l_k}(1)+\sum\limits_{k\in J^-}b_{jk}\omega_{k,l_k}(0),\quad j\in J^-,
\ed
which corresponds to the old boundary condition, and
\be\label{contcond}
\begin{split}
\nu_{j,i}(0)=\nu_{j,i-1}(1), \quad j\in J^+,i=2,\ldots,l_j,\\
\omega_{j,i}(1)=\omega_{j,i-1}(0),\quad j\in J^-,i=2,\ldots,l_j,
\end{split}
\ee
which describes the continuity condition in the artificial vertices. By Theorem \ref{genthm}, to each of the problems \eqref{syst1}, \eqref{ibvpunitvel} there correspond  $C_0$-semigroups $(G(t))_{t\geq 0},$ generated by the operator $A$, and $(\mc G(\tau))_{\tau\geq 0,}$ generated by the operator $\mc A$, where $A$ and $\mc A$ are defined as in \eqref{generatordefinition}. Since we do not distinguish between different values of $p$, we dropped the index $p$ for clarity of notation. The operator $\mb\Q$ provides a similarity relation between these semigroups.
\begin{theorem}\label{simrelthm}
The operator $\mb\Q$ is an isomorphism such that
\be\label{simrel}
G(t)=\mb\Q^{-1}\mc G(ct)\mb\Q,\quad t\geq 0.
\ee
\end{theorem}
\begin{proof}
It is easy to see that $\mb\Q$ is a bounded linear bijection, with the inverse operator given by
\bd
\left(\mb\Q^{-1}\binom{\mb\nu}{\mb\omega}\right)_j(x)=\nu_{j,i}(cL_j(x)-i+1)
\ed
for $x\in[L_j^{-1}((i-1)\slash c),L_j^{-1}(i\slash c)], j\in J^+, i\in\{1,\ldots,l_j\}$, and
\bd
\left(\mb\Q^{-1}\binom{\mb\nu}{\mb\omega}\right)_j(x)=\omega_{j,i}(cL_j(x)-l_j+i)
\ed
for $x\in[L_j^{-1}((l_j-i)\slash c),L_j^{-1}((l_j-i+1)\slash c)], j\in J^-, i\in\{1,\ldots,l_j\}$. The proof of the similarity relation is divided into three steps:
\begin{itemize}
    \item[(i)]  $\mb\Q(D(A))\subset D(\mc A)$ and $\mb\Q^{-1}(D(\mc A))\subset D(A)$,
    \item[(ii)] the operator $\mb\Q^{-1}c\mc A\mb\Q$ is the generator of the semigroup $(\mb\Q^{-1}\mc G(ct)\mb\Q)_{t\geq 0}$,
    \item[(iii)] $\mb\Q^{-1}c\mc A\mb\Q=A$.
\end{itemize}
Let us proceed with (i). Take $(\mb\ups,\mb\w)^T\in D(A)$ and fix $1\leq j\leq 2m$. Then, $\upsilon_j \in W^p_1(0,1)$ by definition and hence \cite[Proposition 9.6]{Brez} implies $\nu_{j,i}\in W^p_1(0,1)$ for any $1\leq i\leq l_j$ due to the fact that $L_j$ is a diffeomorphism. The continuity condition \eqref{contcond} is straightforward. Hence, $\mb\Q(\mb\ups,\mb\w)^T\in D(\mc A)$. The second inclusion in (i) is proved similarly.

Item (ii) is an immediate consequence of \cite[II.2.1 \& II.2.2]{EN}.

It remains to show the equality of the generators. We begin with the domains. Since applying $\mb\Q^{-1}c$ does not affect the domain of the operator $\mb\Q^{-1}c\mc A\mb\Q$, it suffices to consider $\mc A\mb\Q$. Take any function $(\mb\ups,\mb\w)^T\in D(A)$. By (i), $\mb\Q(\mb\ups,\mb\w)^T\in D(\mc A)$, hence $(\mb\ups,\mb\w)^T\in D(\mc A\mb\Q)$ and $D(A)\subset D(\mc A\mb\Q)$. Similarly, if $(\mb\nu,\mb\omega)^T\in D(\mc A\mb\Q)$, then, in particular, $\mb\Q(\mb\nu,\mb\omega)^T\in D(\mc A)$ and by (i), $\mb\Q^{-1}(\mb\Q(\mb\nu,\mb\omega)^T)=(\mb\nu,\mb\omega)^T\in D(A)$ implying $D(\mc A\mb\Q)\subset D(A)$, which proves the equality of domains. Now, take $(\mb\ups,\mb\w)^T\in D(A), j\in J^+$ and $i\in\{1,\ldots,l_j\}$. We have
\begin{displaymath}
\begin{gathered}
\left(c\mc A\mb\Q\binom{\mb\ups}{\mb\w}\right)_{j,i}(y)=c\frac{d}{dy}\upsilon_j\left(L_j^{-1}\left(\frac{y+i-1}{c}\right)\right)\\
=c_j\left(L_j^{-1}\left(\frac{y+i-1}{c}\right)\right)\upsilon_j'\left(L_j^{-1}\left(\frac{y+i-1}{c}\right)\right).
\end{gathered}
\end{displaymath}
Applying the operator $\mb\Q^{-1}$ yields
\bd
\left(\mb\Q^{-1}c\mc A\mb\Q\binom{\mb\ups}{\mb\w}\right)_j(x)=c_j(x)\upsilon_j'(x)=\left(A\binom{\mb\ups}{\mb \w}\right)_j(x)
\ed
for $x\in[L_j^{-1}((i-1)\slash c),L_j^{-1}(i\slash c)]$. Similar calculations for $j\in J^-$ finish the proof of (iii).\\
The similarity relation \eqref{simrel} is a straightforward consequence of the equality of the generators.
\end{proof}

Due to Theorem \ref{simrelthm}, for the remaining part of the paper we assume
\bd
\forall_{j\in\{1,\ldots,2m\}}\quad c_j\equiv 1.
\ed

\subsection{An explicit formula}\label{expforI}
In this subsection we construct an explicit formula for the semigroup governing the solution to the principal part of \eqref{syst1}, that is, with $\mc K=0$. Let us start with necessary notation. If $\mb f:[0,1]\to\R^{2m}$, then we write
\bd
\mb f=\binom{\mb f^+}{\mb f^-},
\ed
where $\mb f^+:[0,1]\to\R^{|J^+|}$ and $\mb f^-:[0,1]\to\R^{|J^-|}$. Similarly, for any $(|J^+|+|J^-|)\times(|J^+|+|J^-|)$ matrix $M$ we write
\bd
M=\left(\begin{array}{cc} M^{11}&M^{12}\\M^{21}&M^{22}\end{array}\right),
\ed
where $M^{11},M^{12},M^{21},M^{22}$ are submatrices of the dimension $|J^+|\times|J^+|,|J^+|\times|J^-|,|J^-|\times|J^+|,|J^-|\times|J^-|$, respectively. These blocks of the boundary matrix $\mc B$ play an essential role in the explicit formula for the semigroup. Further, let us denote
\bd
\mr{\mb\Phi}(x)=\binom{\mr{\mb\ups}(x)}{\mr{\mb\w}(x)},\quad 0<x<1.
\ed
 The formula for the semigroup is found using the method of characteristics, hence we recall here the formulae for the solutions to scalar transport problems. The solution to
\bd
\partial_t \ups+\partial_x \ups=0, \quad \ups(x,0)=\mathring{\ups}(x), \quad \ups(0,t)=\varphi(t),
\ed
$t>0, 0<x<1,$ is given by
\be\label{scalartoright}
\ups(x,t)=\left\{\begin{array}{ll} \mr \ups(x-t),&0< t< x,\\\varphi(t-x),&x< t< x+1,\end{array}\right.
\ee
and, similarly, to
\bd
\partial_t \w-\partial_x \w=0, \quad \w(x,0)=\mathring{\w}(x), \quad \w(1,t)=\psi(t),
\ed
 by
\be\label{scalartoleft}
\w(x,t)=\left\{\begin{array}{ll} \mr \w(x+t),&0< t< 1-x,\\\psi(x+t-1),&1-x< t< 2-x.\end{array}\right.
\ee
Let us return to the problem \eqref{syst1}. For small times, \eqref{syst1} describes a decoupled transport process and hence the solution is given by the translation of the initial condition, that is,
\bd
\binom{\mb\ups^{(0)}(x,t)}{\mb\w^{(0)}(x,t)}=\binom{\mr{\mb\ups}(x-t)}{\mr{\mb\w}(x+t)},
\ed
where $\mb\ups^{(0)}(x,t)$ is defined for $0<t<x$ and $\mb\w^{(0)}(x,t)$ for $0<t<1-x$. In particular, the values $\mb\ups^{(0)}(1,t)=\mr{\mb\ups}(1-t)$ and $\mb\w^{(0)}(0,t)=\mr{\mb\w}(t)$ are well-defined for $0<t<1$. Applying the boundary condition, we have
\bd
\binom{\mb\ups^{(1)}(0,t)}{\mb\w^{(1)}(1,t)}=\mc B\binom{\mb\ups^{(0)}(1,t)}{\mb\w^{(0)}(0,t)}=\mc B\binom{\mr{\mb\ups}(1-t)}{\mr{\mb\w}(t)}.
\ed
Using \eqref{scalartoright} and \eqref{scalartoleft}, we get
\bd
\binom{\mb\ups^{(1)}(x,t)}{\mb\w^{(1)}(x,t)}=\binom{\mathcal{B}^{11}\mathring{\mb\ups}(1-t+x)+\mathcal{B}^{12}\mathring{\mb\w}(t-x)}{\mathcal{B}^{21}\mathring{\mb\ups}(2-t-x)+\mathcal{B}^{22}\mathring{\mb\w}(t+x-1)}.
\ed
Here, $\mb\ups^{(1)}(x,t)$ is defined for $x<t<x+1$, while $\mb\w^{(1)}(x,t)$ for $1-x<t<2-x$.  Continuing this procedure inductively, we can define a family of operators $(T(t))_{t\geq 0}$ by the formula
\begin{equation}\label{Tt}
\begin{split}
T(t)\binom{\mr{\mb\ups}}{\mr{\mb\w}}(x)&=\binom{(T(t)\mr{\mb\Phi})^+(x)}{(T(t)\mr{\mb\Phi})^-(x)}\\
&=\binom{(\mathcal{B}^n)^{11}\mathring{\mb\ups}\left(n-t+x\right)+(\mathcal{B}^n)^{12}\mathring{\mb\w}\left(1-n+t-x\right)}{(\mathcal{B}^n)^{21}\mathring{\mb\ups}\left(1-t+n-x\right)+(\mathcal{B}^n)^{22}\mathring{\mb\w}\left(t-n+x\right)},
\end{split}\end{equation}
where $(T(t)\mr{\mb\Phi})^+(x)$ is defined for $n-1+x<t<n+x$, and $(T(t)\mr{\mb\Phi})^-$ is defined for $n-x<t<n+1-x$, $n\in\N_0$. Then, we have the following
\begin{theorem}\label{formulatheorem}
The family $T=(T(t))_{t\geq 0}$ is a $C_0$-semigroup whose generator is the operator $(A,D(A))$, that is,
\bd
T(t)=G(t),\quad t\geq 0.
\ed
\end{theorem}
\begin{proof}
Since $(A,D(A))$ is the generator of a $C_0$-semigroup, it suffices to show that the Laplace transform $\mc L\{T\}$ of $(T(t))_{t\geq 0}$ equals the resolvent operator $R(\la,A).$  From the proof of \cite[Theorems 4.1 \& 4.2]{JBAB1}, we know that
\begin{displaymath}
\begin{gathered}
R(\lambda,A)\binom{\mb f^+}{\mb f^-}(x)=\left(\begin{array}{cc} \mb{\msf e}_+^{-\la x} & 0\\ 0 & \mb{\msf e}_-^{\la(x-1)}\end{array}\right)\sum\limits_{n=0}^{\infty}\left(\mathcal{B}\mb{\msf e}^{-\la}\right)^n\mathcal{B}\left(\begin{array}{c} \int\limits_0^1\mb{\msf e}_+^{\la(s-1)}\mb f^+(s)\,ds \\ \int\limits_0^1\mb{\msf e}_-^{-\la s}\mb f^-(s)\,ds \end{array}\right)\\
+\left(\begin{array}{c} \int\limits_0^x\mb{\msf e}_+^{\la(s-x)}\mb f^+(s)\,ds \\ \int\limits_x^1\mb{\msf e}_-^{\la(x-s)}\mb f^-(s)\,ds \end{array}\right),
\end{gathered}
\end{displaymath}
where we defined
\bd
\mb{\msf e}_+^z:=\operatorname{diag}(e^z)_{j\in J^+},\quad \mb{\msf e}_-^z:=\operatorname{diag}(e^z)_{j\in J^-}, \quad \mb{\msf e}^z:=\operatorname{diag}(\mb{\msf e}_+^z,\mb{\msf e}_-^z).
\ed
By \cite[Theorem 4.2]{JBAB1}, it suffices to prove the result in $\mb X_1$.

We have, by \eqref{Tt},
\begin{displaymath}
\begin{gathered}
\mc L\{T\}(\lambda)=\int\limits_0^{\infty}\mb{\msf e}^{-\lambda t}T(t)\binom{\mb f^+}{\mb f^-}\,dt
=\left(\begin{array}{cc}\int\limits_0^{1}\mb{\msf e}_+^{-\la(x-s)}\mb f^+(s)\,ds \\ \int\limits_0^{1}\mb{\msf e}_+^{-\la(s-x)}\mb f^-(s)\,ds \end{array}\right)\\
+\sum\limits_{n=1}^{\infty}\left(\begin{array}{cc} \int\limits_{0}^{1}\mb{\msf e}_+^{-\la(n-s+x)}(\mathcal{B}^n)^{11}\mb f^+\left(s\right)\,ds+\int\limits_{0}^{1}\mb{\msf e}_+^{-\la(s+x+n-1)}(\mathcal{B}^n)^{12}\mb f^-\left(s\right)\,ds \\ \int\limits_{0}^{1}\mb{\msf e}_-^{-\la(n+1-s-x)}(\mathcal{B}^n)^{21}\mb f^+\left(s\right)\,ds+\int\limits_{0}^{1}\mb{\msf e}_-^{-\la(s-x+n)}(\mathcal{B}^n)^{22}\mb f^-\left(s\right)\,ds \end{array}\right).
\end{gathered}
\end{displaymath}
Changing the summation parameter from $n$ to $n-1$, the second summand can be written as
\begin{displaymath}
\begin{gathered}
\sum\limits_{n=0}^{\infty}\left(\begin{array}{cc} \int\limits_{0}^{1}\mb{\msf e}_+^{-\la(n+1-s+x)}(\mathcal{B}^{n+1})^{11}\mb f^+\left(s\right)\,ds+\int\limits_{0}^{1}\mb{\msf e}_+^{-\la(s+x+n)}(\mathcal{B}^{n+1})^{12}\mb f^-\left(s\right)\,ds \\ \int\limits_{0}^{1}\mb{\msf e}_-^{-\la(n+2-s-x)}(\mathcal{B}^{n+1})^{21}\mb f^+\left(s\right)\,ds+\int\limits_{0}^{1}\mb{\msf e}_-^{-\la(s-x+n+1)}(\mathcal{B}^{n+1})^{22}\mb f^-\left(s\right)\,ds \end{array}\right)\\
=\sum\limits_{n=0}^{\infty}\left(\begin{array}{cc} \mb{\msf e}_+^{-\la x}\mb{\msf e}_+^{-\la n}(\mathcal{B}^{n+1})^{11} & \mb{\msf e}_+^{-\la x}\mb{\msf e}_+^{-\la n}(\mathcal{B}^{n+1})^{12} \\ \mb{\msf e}_-^{-\la (1-x)}\mb{\msf e}_-^{-\la n}(\mathcal{B}^{n+1})^{21} & \mb{\msf e}_-^{-\la (1-x)}\mb{\msf e}_-^{-\la n}(\mathcal{B}^{n+1})^{22} \end{array}\right)\left(\begin{array}{c} \int\limits_0^1\mb{\msf e}_+^{-\la (1-s)}\mb f^+\left(s\right)\,ds \\ \int\limits_0^1\mb{\msf e}_-^{-\la s}\mb f^-\left(s\right)\,ds \end{array}\right)\\
=\sum\limits_{n=0}^{\infty}\left(\begin{array}{cc} \mb{\msf e}_+^{-\la x}\mb{\msf e}_+^{-\la n} & 0 \\ 0 & \mb{\msf e}_-^{-\la (1-x)}\mb{\msf e}_-^{-\la n} \end{array}\right)\mathcal{B}^{n+1}\left(\begin{array}{c} \int\limits_0^1\mb{\msf e}_+^{\la (s-1)}\mb f^+\left(s\right)\,ds \\ \int\limits_0^1\mb{\msf e}_-^{-\la s}\mb f^-\left(s\right)\,ds \end{array}\right)\\
=\left(\begin{array}{cc} \mb{\msf e}_+^{-\la x} & 0 \\ 0 & \mb{\msf e}_-^{\la (x-1)} \end{array}\right)\sum\limits_{n=0}^{\infty}\mb{\msf e}^{-\la n}\mathcal{B}^{n+1}\left(\begin{array}{c} \int\limits_0^1\mb{\msf e}_+^{\la (s-1)}\mb f^+\left(s\right)\,ds \\ \int\limits_0^1\mb{\msf e}_-^{-\la s}\mb f^-\left(s\right)\,ds \end{array}\right)\\
=\left(\begin{array}{cc} \mb{\msf e}_+^{-\la x} & 0 \\ 0 & \mb{\msf e}_-^{\la (x-1)} \end{array}\right)\sum\limits_{n=0}^{\infty}\left(\mathcal{B}\mb{\msf e}^{-\la }\right)^n\mathcal{B}\left(\begin{array}{c} \int\limits_0^1\mb{\msf e}_+^{\la (s-1)}\mb f^+\left(s\right)\,ds\\\int\limits_0^1\mb{\msf e}_-^{-\la s}\mb f^-\left(s\right)\,ds\end{array}\right),
\end{gathered}
\end{displaymath}
since the matrix $\mb{\msf e}^{-\la n}$ commutes with $\mc B$. Comparing the formulae, we obtain $\mc L\{T\}(\lambda) = R(\la, A)$.
\end{proof}

\subsection{Reduction to the transport semigroup}
If we take $J^-=\emptyset$, then \eqref{syst1} becomes a pure transport problem, which was extensively investigated in \cite{KS2005,Dorn,BN} in the case $p=1$. In \cite[Prop. 3.3]{Dorn} it was shown that if $\mc B$ is stochastic, then the family $(S(t))_{t\geq 0}$ defined by
\be\label{transportsemigroup}
S(t)\binom{\mb\ups}{\mb\w}(x)=\mc B^n\binom{\mb\ups}{\mb\w}(n-t+x),\quad 0<n-t+x<1, n\in\N_0,
\ee
is a $C_0$-semigroup on $\mb X_1$, which governs the solution to the transport problem. It is not difficult to show, \cite[Theorem 3.1]{BFN3}, that it is also a $C_0$-semigroup for an arbitrary matrix $\mc B$ as well. Since $(S(t))_{t\geq 0}$ is a special case of $(G(t))_{t\geq 0}$, the family $(S(t))_{t\geq 0}$ is also the transport semigroup on any $\mb X_p$, $1\leq p<\infty$. We shall show that $(G(t))_{t\geq 0}$ is similar to $(S(t))_{t\geq 0}$. The advantage of this similarity relation is that, by allowing for a straightforward application of the spectral decomposition of the matrix $\mc B$, it significantly simplifies the analysis of the long term asymptotics of $(G(t))_{t\geq 0}$ given in Section \ref{asymptoticbehaviour}.

Define $\mb{\mb{\mc V}}:\mb X_p\to\mb X_p$ by
\bd
\mb{\mc V}\binom{\mb\ups}{\mb\w}(x)=\binom{\mb\ups(x)}{\mb\w(1-x)}.
\ed
The action of this operator  reverses the direction of the edges with $j\in J^-$ so that the flows occur now from $0$ to $1$ on all edges.
\begin{proposition}\label{sgsim}
The operator $\mb{\mc V}$ is an isometric isomorphism in any $\mb X_p$ with
\begin{equation}
\mb{\mc V}^{-1} = \mb{\mc V},
\label{vinv}
\end{equation}
satisfying
\begin{equation}\label{simrel1}
G(t)=\mb{\mc V}S(t)\mb{\mc V},\quad t\geq 0.
\end{equation}
In particular, for $0< n-t+x<1, n\in\N_0, $
\begin{equation}
G(t)\binom{\mb\ups}{\mb\w}(x)=\mb{\mc V}\mc B^n\mb{\mc V}\binom{\mb\ups}{\mb\w}(n-t+x).
\label{GS}
\end{equation}
\end{proposition}
\begin{proof}
That $\mb{\mc V}$ is an isomorphism satisfying \eqref{vinv} is clear. For \eqref{simrel1}, we have
\begin{displaymath}
\begin{gathered}
\mb{\mc V}G(t)\binom{\mb\ups}{\mb\w}(x)=\binom{(\mb{\mc V}G(t)\mb\Phi)^+(x)}{(\mb{\mc V}G(t)\mb\Phi)^-(x)}\\
=\binom{(\mathcal{B}^n)^{11}\mb\ups\left(n-t+x\right)+(\mathcal{B}^n)^{12}\mb\w\left(t-n+1-x\right)}{(\mathcal{B}^n)^{21}\mb\ups\left(n-t+x\right)+(\mathcal{B}^n)^{22}\mb\w\left(t-n+1-x\right)}.
\end{gathered}
\end{displaymath}
Observe that the formula for $(\mb{\mc V}G(t)\mb\Phi)^+(x)$ is valid for $n-1+x<t<n+x$, while $(\mb{\mc V}G(t)\mb\Phi)^-(x)$ is defined for $n-(1-x)<t<n+1-(1-x)$. Both of these conditions yield $0<n-t+x<1,$ as in \eqref{transportsemigroup}. Hence, we can write
\begin{equation}
\mb{\mc V}G(t)\binom{\mb\ups}{\mb\w}(x)=\mc B^n\mb{\mc V}\binom{\mb\ups}{\mb\w}(n-t+x)=S(t)\mb{\mc V}\binom{\mb\ups}{\mb\w}(x).
\label{GS1}
\end{equation}
\end{proof}

\section{Asymptotic behaviour}\label{asymptoticbehaviour}
Let $\sigma({\mc B})=\{\la_1,\ldots,\la_k\}, k\leq 2m$, be the set of the eigenvalues of the matrix $\mc B$. For any $i\in\{1,\ldots,k\},$ denote by $\alpha_i$ the algebraic multiplicity of $\la_i$. Further, let $\{\mb E_{i_1},\ldots,\mb E_{i_{\alpha_i}}\}$ and $\{\mb F_{i_1},\ldots,\mb F_{i_{\alpha_i}}\}$ be the sets of right and left (generalized) eigenvectors corresponding to the eigenvalue $\la_i$, respectively, selected so as $\mb F_{i_j}\cdot \mb E_{i_l}=\delta_{jl}$ for any $j,l\in\{1,\ldots,\alpha_i\}$. Then, denoting by $\Pi_i$ the spectral projection onto the right eigenspace $Lin\{\mb E_{i_1},\ldots,\mb E_{i_{\alpha_i}}\}$, for any vector $\mb U\in\R^{2m}$ we have
\begin{equation}
\mb U=\sum\limits_{i=1}^k\Pi_i\mb U=\sum\limits_{i=1}^k\left(\sum\limits_{j=1}^{\alpha_i}\left(\mb F_{i_j}\cdot\mb U\right)\mb E_{i_j}\right).
\label{Urep}
\end{equation}
Since the matrices $\Pi_j, 1\leq j\leq k$, form the spectral resolution of the identity, that is, $\sum_{j=1}^k\Pi_j=I$, using the binomial expansion, for $n\in\N_0$ we have
\begin{equation}\label{bpowers}
\begin{split}
\mc B^n=\sum\limits_{j=1}^k\mc B^n\Pi_j=\sum\limits_{j=1}^k(\lam_jI+\mc B-\lam_jI)^n\Pi_j\\
=\sum\limits_{j=1}^k\sum\limits_{r=0}^n\binom{n}{r}\lam_j^{n-r}(\mc B-\lam_jI)^r\Pi_j=\sum\limits_{j=1}^k\lam_j^n\mc p_j(n)\Pi_j,
\end{split}
\end{equation}
where $\mc p_j$ is a matrix-valued polynomial in $n$ of the degree strictly smaller than $\alpha_j$.

 Let us take eigenvalues $\lam_{j_1},\ldots,\lam_{j_h}\in\sigma(\mc B)$ for some  $1\leq h\leq k$ and  consider a family $(G^{\lam}(t))_{t\geq 0}$ of linear and bounded operators given by
\begin{equation}
G^{\lam}(t)\binom{\mb\ups}{\mb \w}(x):=\mb{\mc V}\sum\limits_{r=1}^h\lam_{j_r}^n\mc p_{j_r}(n)\Pi_{j_r}\mb{\mc V}\binom{\mb\ups}{\mb \w}(n-t+x),
\label{glam}
\end{equation}
defined for $0< n-t+x< 1, n\in\N_0$.

\subsection{Invariant subspaces}
First,  we show that the family $(G^{\lam}(t))_{t\geq 0}$ has an invariant subspace. We begin with the following lemma.
\begin{lemma}
The family $(G^{\lam}(t))_{t\geq 0}$ is a semigroup.
\end{lemma}
\begin{proof}
We see that \eqref{glam} can be written, by \eqref{transportsemigroup} and \eqref{bpowers}, as
\begin{equation}
G^{\lam}(t)=\sum\limits_{r=1}^h\mb{\mc V}S(t)\Pi_{j_r}\mb{\mc V},
\label{glam1}
\end{equation}
where we identify the matrix multiplication by $\Pi_{j_r}$ with the operator on $\mb X_p$ defined by $[\Pi_{j_r}\mb u](x) = \Pi_{j_r}\mb u(x), \mb u \in \mb X_p$, a.e. $x\in [0,1]$. Now, using again  \eqref{transportsemigroup} and \eqref{bpowers}, $\Pi_{j}S(t)= S(t)\Pi_{j}$ and hence $S(t+s)\Pi_{j} = S(t)S(s)\Pi^2_{j} = S(t)\Pi_{j}S(s)\Pi_{j}$,  for each $1\leq j\leq k$, $t,s\geq 0$. Thus, for any $1\leq i,j\leq k$,
$$
\mb{\mc V}S(t)\Pi_{j} \mb{\mc V}\mb{\mc V} S(s)\Pi_{i}\mb{\mc V} = \mb{\mc V}S(t)\Pi_{j}\Pi_i S(s)\mb{\mc V}= \left\{\begin{array}{lcl}\mb{\mc V}S(t+s)\Pi_{j}\mb{\mc V}&\text{if}&i=j,\\
0&\text{if}&i\neq j.\end{array}\right.
$$
Then,
$$
G^{\lam}(t)G^{\lam}(s) = \sum\limits_{r=1}^h\mb{\mc V}S(t+s)\Pi_{j_r}\mb{\mc V} = G^{\lam}(t+s).
$$
\end{proof}
The last formula implies, in particular, that
\begin{equation}\label{sgprop0}
    G^{\lam}(0)G^{\lam}(t)=G^{\lam}(t)G^{\lam}(0)=G^{\lam}(t), \quad t\geq 0.
\end{equation}
In other words, for any $t\geq 0$ and $(\mb\ups, \mb\w)^T\in \mb X_p$,
$$
G^{\lam}(t)\binom{\mb\ups}{\mb\w}\in \operatorname{rng}G^{\lam}(0),
$$
where "rng" denotes the range of an operator. This shows that the range of $G^{\lam}(0)$ is a candidate for an invariant subspace for the semigroup $(G(t))_{t\geq 0}$. We shall investigate this observation further. Let us define sets
$$
Z_j:=\left\{\mb{\mc V}\Pi_j\mb{\mc V}\binom{\mb\ups}{\mb\w}:\binom{\mb\ups}{\mb\w}\in\mb X_p\right\},\quad j=1,\ldots,k.
$$
Certainly, each  set is a linear subspace of the space $\mb X_p.$  Moreover, $Z_{j_1}\cap Z_{j_2}=\{0\}$ for $j_1\neq j_2$. Indeed, take any $(\mb\ups,\mb\w)^T\in Z_{j_1}\cap Z_{j_2}$. Then there exist $(\mb\ups_1,\mb\w_1)^T\in Z_{j_1}, (\mb\ups_2,\mb\w_2)^T\in Z_{j_2}$ such that
$$
\mb{\mc V}\Pi_{j_1}\mb{\mc V}\binom{\mb\ups_1}{\mb\w_1}=\binom{\mb\ups}{\mb\w}=\mb{\mc V}\Pi_{j_2}\mb{\mc V}\binom{\mb\ups_2}{\mb\w_2}.
$$
Applying the operator $\mb{\mc V}$ and multiplying by $\Pi_{j_1}$ we obtain $$
0=\Pi_{j_1}^2\mb{\mc V}\binom{\mb\ups_1}{\mb\w_1}=\Pi_{j_1}\mb{\mc V}\binom{\mb\ups_1}{\mb\w_1} = \mb{\mc V}\binom{\mb\ups}{\mb\w}.
$$
Since $\mb{\mc V}$ is an isomorphism, $(\mb\ups,\mb\w)^T=0$.

By the definition of the image of a map we have
\begin{displaymath}
\begin{gathered}
Z^{\lam}:=\operatorname{rng}G^{\lam}(0)=\left\{\mb{\mc V}\sum\limits_{r=1}^h\Pi_{j_r}\mb{\mc V}\binom{\mb\ups}{\mb\w}:\binom{\mb\ups}{\mb\w}\in\mb X_p\right\}=\bigoplus\limits_{r=1}^hZ_{j_r},
\end{gathered}
\end{displaymath}
where the fact that the sum is direct follows as above. We observe also that the subspace $Z^{\lam}$ is closed since the operators induced by the matrices $\Pi_{j_r}$ are also continuous projectors on $\mb X_p$ and $\mb{\mc V}$ is an isomorphism.
\begin{lemma}
The subspace $Z^{\lam}$ is invariant under $(G(t))_{t\geq 0}$.
\end{lemma}
\begin{proof}
Take any element from the space $Z^{\lam}$, which is of the form $\mb{\mc V}\sum_{r=1}^h\Pi_{j_r}\mb{\mc V}(\mb\ups,\mb\w)^T$ for some $(\mb\ups,\mb\w)^T\in\mb X_p$. We have, by \eqref{GS1}, \eqref{glam1} and \eqref{sgprop0},
\begin{displaymath}
\begin{gathered}
G(t)\left(\mb{\mc V}\sum\limits_{r=1}^h\Pi_{j_r}\mb{\mc V}\binom{\mb\ups}{\mb \w}\right)=\mb{\mc V}S(t)\mb{\mc V}\left(\mb{\mc V}\sum\limits_{r=1}^h\Pi_{j_r}\mb{\mc V}\binom{\mb\ups}{\mb \w}\right)\\
=G^{\lam}(t)\binom{\mb\ups}{\mb \w}
=G^{\lam}(0)G^{\lam}(t)\binom{\mb\ups}{\mb \w}\in\operatorname{rng}G^{\lam}(0)=Z^{\lam}.
\end{gathered}
\end{displaymath}
\end{proof}
From the above proof we see also that the family $(G^{\lam}(t))_{t\geq 0}$ is the restriction of the semigroup $(G(t))_{t\geq 0}$ to the subspace $Z^{\lam}$. Since this subspace is invariant and closed, we have proved the following
\begin{proposition}\label{subspacesemigroup}
The family $(G^{\lam}(t))_{t\geq 0}$ is a $C_0$-semigroup on the space $Z^{\lam}=\operatorname{rng}G^{\lam}(0)$.
\end{proposition}

\subsection{Asymptotic stability}
Let us assume that $\sigma(\mc B)$ is ordered as $\sigma(\mc B)=\{\lam_1,\ldots,\lam_d,\lam_{d+1},\ldots,\lam_k\}$. We consider a family $(G^2(t))_{t\geq 0}$ given by
\begin{displaymath}
G^2(t)\binom{\mb\ups}{\mb\w}(x)=\mb{\mc V}\sum\limits_{j=d+1}^k\lam_j^n\mc p_j(n)\Pi_j\mb{\mc V}\binom{\mb\ups}{\mb\w}(n-t+x),
\end{displaymath}
defined for $0< n-t+x<1, n\in\N_0$. By Proposition \ref{subspacesemigroup}, $(G^2(t))_{t\geq 0}$ is a $C_0$-semigroup on the space $Z^2:=\operatorname{rng}G^2(0)$.
\begin{lemma}\label{aslemma}
If \begin{equation}\label{stableeigenvalues}
|\lam_j|<1,\quad j=d+1,\ldots,k,
\end{equation}
then the semigroup $(G^2(t))_{t\geq 0}$ is uniformly asymptotically stable.
\end{lemma}
The proof is similar to that of \cite[Thm. 5.1(iii)]{BN}, which, although carried out for the $L^1$ case, can be easily adapted to any $L^p$ space. Before we proceed with the proof, let us clarify that in the space $\R^{2m}$ we consider the $p$-norm $|x|=(\sum_{i=1}^{2m}|x_i|^p)^{\frac{1}{p}}$, and for any $2m\times 2m$ matrix $A$ by $|A|$ we denote the norm of the linear operator induced by $A$ corresponding to the $p$-norm. This notation does not lead to any confusion since the distinction between a vector and a matrix is clear from the context. In particular, for any $\mb f\in\mb X_p$,
$$
\left\|\mb f\right\|=\left(\int\limits_0^1|\mb f(x)|^p\,dx\right)^{\frac{1}{p}}=\left(\int\limits_0^1\sum\limits_{i=1}^{2m}|f_i(x)|^p\,dx\right)^{\frac{1}{p}}.
$$
\begin{proof}
First, let us take $j\in\{d+1,\ldots,k\}$ and let $\tau_n=t-n$ be such that $0\leq \tau_n\leq 1$. Using the fact that $|\lambda_j|<1,$ we can find $0<\bar\la_j<1$ such that ${|\la_j|}/{\bar\la_j}<1$. Taking $\la=\max_j {|\la_j|}/{\bar\la_j}$ we have that ${|\la_j|}/{\bar\la_j}\leq\la<1$ for all $d+1\leq j\leq k.$ Since $|\bar\la^n\mc p_j(n)\Pi_j|$ is uniformly bounded with respect to $n\in\N_0$, denoting $\theta(n,x,t) = n-t+x$, we have
\begin{align*}
&\left\|\lambda_j^n\mc p_j(n)\Pi_j\mb{\mc V}\binom{\mb\ups}{\mb\w}(n-t+\cdot)\right\|^p \leq \lambda^{n}\int\limits_{\tau_n}^{1}\left|\bar\la_j^n\mc p_j(n)\Pi_j\mb{\mc V}\binom{\mb\ups}{\mb\w}(\theta(n,x,t))\right|^p\,dx\\&\phantom{xxxxxxxxxxxxxxxxxx}+\lambda^{n+1}\int\limits_{0}^{\tau_n}\left|\bar\la_j^{n+1}\mc p_j(n+1)\Pi_j\mb{\mc V}\binom{\mb\ups}{\mb\w}(\theta(n+1,x,t))\right|^p\,dx\\
&\leq C_{\mc p_j}^p\lambda^{np}\left(\!\int\limits_0^{\tau_n}\left|\binom{\mb\ups(\theta(n+1,x,t))}{\mb\w(1-\theta(n+1,x,t))}\right|^p\!\!\!dx
+\int\limits_{\tau_n}^1\left|\binom{\mb\ups(\theta(n,x,t))}{\mb\w(1-\theta(n,x,t))}\right|^p\!\!\!dx\!\right)\\
&=C_{\mc p_j}^p\lambda^{np}\left(\int\limits_{-\tau_n+1}^1\left|\binom{\mb\ups(s)}{\mb\w(1-s)}\right|^p\,ds+
\int\limits_0^{-\tau_n+1}\left|\binom{\mb\ups(s)}{\mb\w(1-s)}\right|^p\,ds\right)\\
&=C_{\mc p_j}^p\lambda^{np}\left\|\binom{\mb\ups}{\mb\w}\right\|^p.
\end{align*}
Moreover, since  $0< n-t+x< 1$ implies $n\leq t+1$,  we can write
\begin{displaymath}
\begin{gathered}
\left\|G^2(t)\binom{\mb\ups}{\mb \w}\right\|\leq\left(\sum\limits_{j=d+1}^kC_{\mc p_j}\right)\lambda^n\left\|\binom{\mb\ups}{\mb\w}\right\|\leq\left(\sum\limits_{j=d+1}^kC_{\mc p_j}\lambda\right)e^{t\ln{\lambda}}\left\|\binom{\mb\ups}{\mb\w}\right\|,
\end{gathered}
\end{displaymath}
with $\ln \la<0$, which shows the uniform asymptotic stability of $(G^2(t))_{t\geq 0}$.
\end{proof}

\subsection{Periodicity}
Similarly to the previous paragraph, we define a family $(G^1(t))_{t\geq 0}$ by
\begin{displaymath}
G^1(t)\binom{\mb\ups}{\mb\w}(x)=\mb{\mc V}\sum\limits_{j=1}^d\lam_j^n\mc p_j(n)\Pi_j\mb{\mc V}\binom{\mb\ups}{\mb\w}(n-t+x),
\end{displaymath}
$0< n-t+x<1, n\in\N_0$. Again, $(G^1(t))_{t\geq 0}$ is a $C_0$-semigroup on the space $Z^1:=\operatorname{rng}G^1(0)$. Let us assume that the eigenvalues $\lam_j, j=1,\ldots,d$, are semisimple of the form
\begin{equation}\label{periodiceigenvalues}
\lam_j=e^{\frac{2\pi i(j-1)}{d}},
\end{equation}
where $i$ is the imaginary unit. This implies that $\mc B\Pi_j=\lam_j\Pi_j, j=1,\ldots,d$.

\begin{lemma}\label{perlemma}
The semigroup $(G^1(t))_{t\geq 0}$ is periodic with period $d$.
\end{lemma}
\begin{proof}
In the first step we prove that the semigroup $(G^1(t))_{t\geq 0}$ is periodic. Put $n'=n+d$. Then, $0< n'-(t+d)+x<1$ and, by \eqref{periodiceigenvalues}, we have
\begin{displaymath}
\begin{gathered}
G^1(t+d)\binom{\mb\ups}{\mb\w}(x)=\mb{\mc V}\sum\limits_{j=1}^d\lam_j^{n'}\Pi_j\mb{\mc V}\binom{\mb\ups}{\mb\w}(n'-t-d+x)\\
=\mb{\mc V}\sum\limits_{j=1}^de^{\frac{2\pi i(j-1)n}{d}}\Pi_j\mb{\mc V}\binom{\mb\ups}{\mb\w}(n-t+x)=G^1(t)\binom{\mb\ups}{\mb\w}(x).
\end{gathered}
\end{displaymath}
Hence, $(G^1(t))_{t\geq 0}$ is indeed periodic and its period does not exceed $d$. To prove that it equals $d$, let us calculate the Laplace transform of $(G^1(t))_{t\geq 0}$. We have
\begin{displaymath}
\begin{gathered}
\int\limits_0^{\infty}e^{-\mu t}G^1(t)\binom{\mb\ups}{\mb\w}(x)\,dt=\int\limits_0^{\infty}e^{-\mu t}\mb{\mc V}\sum\limits_{j=1}^de^{\frac{2\pi i(j-1)n}{d}}\Pi_j\mb{\mc V}\binom{\mb\ups}{\mb\w}(n-t+x)\,dt\\
=\mb{\mc V}\sum\limits_{j=1}^{d}\left(\int\limits_0^xe^{-\mu(x-s)}\Pi_j\binom{\mb\ups(s)}{\mb\w(1-s)}\,ds\right.\\
\left.+\sum\limits_{n=1}^{\infty}e^{\left(\frac{2\pi i(j-1)}{d}-\mu\right)n}\int\limits_{0}^{1}e^{-\mu(x-s)}\Pi_j\binom{\mb\ups(s)}{\mb\w(1-s)}\,ds\right).
\end{gathered}
\end{displaymath}
Since
$$
\sum\limits_{n=1}^{\infty}e^{\left(\frac{2\pi i(j-1)}{d}-\mu\right)n} = \frac{e^{\left(\frac{2\pi i(j-1)}{d}-\mu\right)}}{1-e^{\left(\frac{2\pi i(j-1)}{d}-\mu\right)}} =: f_j(\mu)
$$
for $\Re\mu>0$ and $f_j$ extends to an analytic function on $\mbb C$ except for $\mu  = 2\pi i (j-1)/{d} +2\pi i l$ for any $l\in \mbb Z$ (the set of integers), where it has first order poles,
the resolvent of the generator of  $(G^1(t))_{t\geq 0}$ has singularities only where one of the $f_j$s has a pole. Since any $k\in \mbb Z$ can be written as $k = ld +j-1$ for some $l\in \mbb Z, 1\leq j\leq d$, we see that $2\pi i k/d = 2\pi i l + 2\pi i(j-1)/d$, we see that the resolvent is analytic in $\mbb C$ except for $\mu \in \frac{2\pi i}{d} \cdot\mbb Z$. Hence, by \cite[Lemma IV.2.25]{EN},  we conclude that the period of $(G^1(t))_{t\geq 0}$ equals $d$.
\end{proof}

\subsection{The limit semigroup}
To formulate the main theorem of this section, we note that
$$
G(t)=G^1(t)+G^2(t),\quad t\geq 0,
$$
and, since the matrices $\Pi_j, j=1,\ldots,k$, form the spectral resolution of the identity, there holds
$$
\mb X_p=Z^1\oplus Z^2.
$$
Then, combining Lemmas \ref{aslemma} and \ref{perlemma}, we obtain
\begin{theorem}\label{maintheorem}
Under the assumption \eqref{stableeigenvalues}, there exists a decomposition $\mb X_p=Z^1\oplus Z^2$ into $(G(t))_{t\geq 0}$-invariant subspaces $Z^1$ and $Z^2$ such that $(G(t))_{t\geq 0}$  asymptotically as $t\to\infty$ behaves as  $(G^1(t))_{t\geq 0}=(G(t)|_{Z^1})_{t\geq 0}$ in the sense that
\begin{equation}
\lim\limits_{t\to \infty}(G(t)- G^1(t))= 0
\label{lim1}
\end{equation}
in the uniform operator topology in any $\mb X_p$, $p\in [1,\infty)$. The semigroup $(G^1(t))_{t\geq 0}$ is explicitly given by
\begin{displaymath}
G^1(t)\binom{\mb\ups}{\mb\w}(x)=\binom{(G^1(t)\mb\Phi)^+(x)}{(G^1(t)\mb\Phi)^-(x)},
\end{displaymath}
where
\begin{align*}
(G^1(t)\mb\Phi)^+(x)&=\sum\limits_{j=1}^d\lam_j^n\left((\mc p_j(n)\Pi_j)^{11}\mb\ups(n-t+x)\right.\\
&+\left.(\mc p_j(n)\Pi_j)^{12}\mb\w(1-n+t-x)\right),
\end{align*}
valid for $n-1+x< t<n+x$, and
\begin{align*}
(G^1(t)\mb\Phi)^-(x)&=\sum\limits_{j=1}^d\lam_j^n\left((\mc p_j(n)\Pi_j)^{21}\mb\ups(1-t+n-x)\right.\\
&+\left.(\mc p_j(n)\Pi_j)^{22}\mb\w(t-n+x)\right),
\end{align*}
valid for $n-x< t<n+1-x$. If, in addition, \eqref{periodiceigenvalues} is satisfied, then $(G^1(t))_{t\geq 0}$ is periodic with period $d$ and
\begin{displaymath}
G^1(t)\binom{\mb\ups}{\mb\w}(x)=\left(\begin{array}{cc}\sum\limits_{j=1}^de^{\frac{2\pi i(j-1)n}{d}}\left(\Pi_j^{11}\mb\ups(n-t+x)+\Pi_j^{12}\mb\w(1-n+t-x)\right)\\
\sum\limits_{j=1}^de^{\frac{2\pi i(j-1)n}{d}}\left(\Pi_j^{21}\mb\ups(1-t+n-x)+\Pi_j^{22}\mb\w(t-n+x)\right)
\end{array}\right).
\end{displaymath}
\end{theorem}

\begin{remark}
We note that if $\mc B$ is stochastic, which is the case considered in \cite{BN}, then the assumptions \eqref{stableeigenvalues} and \eqref{periodiceigenvalues} are automatically satisfied by \cite[Chapter 8, p. 696]{Mey}. However, in Section \ref{examplesection} we present an example when the boundary matrix is not stochastic, yet these assumptions are still satisfied.
\end{remark}

\section{Examples}\label{examplesection}

\begin{example}
Let $\mb\ups=(\ups_1,\ups_2), \mb\w=(\w_3,\w_4)$ and consider
\begin{equation}\label{firstexample}
\begin{gathered}
\partial_t\mb\ups(x,t)=-\partial_x\mb\ups(x,t),\quad t>0, 0<x<1,\\
\partial_t\mb\w(x,t)=\partial_x\mb\w(x,t),\quad t>0, 0<x<1,\\
\mb\ups(x,0)=\mr{\mb\ups}(x),\; \mb\w(x,0)=\mr{\mb\w}(x),\quad 0<x<1,\\
\binom{\mb\ups(0,t)}{\mb\w(1,t)}=\mc B\binom{\mb\ups(1,t)}{\mb\w(0,t)},\quad t>0,
\end{gathered}
\end{equation}
on a graph presented in Figure \ref{1stexgraph}, with $l_1(\bv_1)=l_2(\bv_2)=0$ and $l_1(\bv_2)=l_2(\bv_3)=1$.
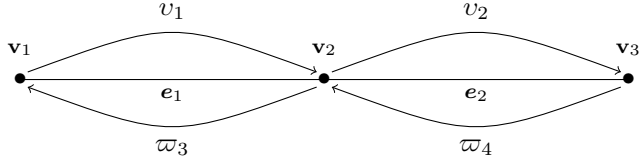
\begin{figure}[h]
\centering
\begin{tikzpicture}
\draw [-] (-4,0) -- (-0,0);
\draw [-] (0.0,0) -- (4.05,0);

\draw (-2,-0.2)node{$\footnotesize{\mb e_1}$};
\draw (2,-0.2)node{$\footnotesize{\mb e_2}$};

\draw (-4.0,0.4) node {\footnotesize{$\bv_1$}};
\draw (-0.0,0.4) node {\footnotesize{$\bv_2$}};
\draw (4.0,0.4) node {\footnotesize{$\bv_3$}};

\draw (0,0)node{$\bullet$};
\draw (4,0)node{$\bullet$};
\draw (-4,0)node{$\bullet$};

\draw (-2,0.9)node{$\footnotesize{\upsilon_1}$};
\draw (-2,-0.9)node{$\footnotesize{\varpi_3}$};
\draw (2,0.9)node{$\footnotesize{\upsilon_2}$};
\draw (2,-0.9)node{$\footnotesize{\varpi_4}$};

\draw[->] (-3.9,0.1)..controls(-2,0.8) ..(-0.1,0.1);
\draw[<-] (-3.9,-0.1)..controls(-2,-0.8) ..(-0.1,-0.1);
\draw[<-] (3.9,0.1)..controls(2,0.8) ..(0.1,0.1);
\draw[->] (3.9,-0.1)..controls(2,-0.8) ..(0.1,-0.1);

\end{tikzpicture}
\caption{The graph for the problem \eqref{firstexample} with the directions of the flows.}
\label{1stexgraph}
\end{figure}
Let
\begin{displaymath}
\mc B=\left(\begin{array}{cccc}
0&0&1&0\\
\frac{1}{4}&0&0&\frac{1}{2}\\
\frac{3}{4}&0&0&\frac{1}{2}\\
0&1&0&0
\end{array}\right).
\end{displaymath}
Then, we have $\lam_1=1,\lam_2=-1,\lam_3=\frac{1}{2},\lam_4=-\frac{1}{2}$ and the right and left eigenvectors corresponding to the peripheral spectrum are $\mb E_1=(2,1,2,1)^T,\mb F_1=\frac{1}{6}(1,1,1,1)^T,\mb E_2=(2,-1,-2,1)^T$ and $\mb F_2=\frac{1}{6}(1,-1,-1,1)^T$. We calculate the projections to be
\begin{displaymath}
\Pi_1=\frac{1}{6}
\left(\begin{array}{cccc}
2&2&2&2\\
1&1&1&1\\
2&2&2&2\\
1&1&1&1
\end{array}\right),\quad
\Pi_2=\frac{1}{6}
\left(\begin{array}{cccc}
2&-2&-2&2\\
-1&1&1&-1\\
-2&2&2&-2\\
1&-1&-1&1
\end{array}\right).
\end{displaymath}
By Theorem \ref{maintheorem}, the semigroup $(G(t))_{t\geq 0}$ governing the solution to \eqref{firstexample} converges asymptotically to the periodic semigroup $(G^1(t))_{t\geq 0}$ given by
\begin{displaymath}
\begin{gathered}
(G^1(t)\mr{\mb\Phi})^+(x)=\frac{1}{6}\left(\left(\begin{array}{cc}2&2\\1&1\end{array}\right)+(-1)^n\left(\begin{array}{cc}2&-2\\-1&1\end{array}\right)\right)\left(\begin{array}{c}\mr{\ups}_1(n-t+x)\\\mr{\ups}_2(n-t+x)\end{array}\right)\\
+\left(\left(\begin{array}{cc}2&2\\1&1\end{array}\right)+(-1)^n\left(\begin{array}{cc}-2&2\\1&-1\end{array}\right)\right)\left(\begin{array}{c}\mr{\w}_3(1-n+t-x)\\\mr{\w}_4(1-n+t-x)\end{array}\right),\\
(G^1(t)\mr{\mb\Phi})^-(x)=\frac{1}{6}\left(\left(\begin{array}{cc}2&2\\1&1\end{array}\right)+(-1)^n\left(\begin{array}{cc}-2&2\\1&-1\end{array}\right)\right)\left(\begin{array}{c}\mr{\ups}_1(1-t+n-x)\\\mr{\ups}_2(1-t+n-x)\end{array}\right)\\
+\left(\left(\begin{array}{cc}2&2\\1&1\end{array}\right)+(-1)^n\left(\begin{array}{cc}2&-2\\-1&1\end{array}\right)\right)\left(\begin{array}{c}\mr{\w}_3(t-n+x)\\\mr{\w}_4(t-n+x)\end{array}\right).
\end{gathered}
\end{displaymath}
\end{example}

\begin{example}\label{beamexample}
Following \cite[Example 7.1.4]{JaZwbook}, let us consider the model of a Timoshenko beam given by
\begin{equation}\label{secondorderbeam}
\begin{split}
    \rho\p^2_t\mc v(x,t)=K(\p^2_x\mc v(x,t)-\p_x\mc u(x,t)),\\
    I_{\rho}\p^2_t\mc u(x,t)=EI\p^2_x\mc u(x,t)+K(\p_x\mc v(x,t)-\mc u(x,t)),
\end{split}
\quad x\in(0,l), t\geq 0,
\end{equation}
where $\mc v$ is the transverse displacement of the beam and $\mc u$ is the rotation angle of a cross-sectional filament of the beam, while the coefficients $\rho, I_{\rho}, EI$ and $K$ are the mass per unit length, the rotary moment of inertia of a cross section, the product of Young's modulus of elasticity and the moment of inertia of a cross section, and the shear modulus, respectively. As proposed in \cite{JaZwbook}, introducing
\begin{displaymath}
\left.\begin{array}{lll}
p_1=\p_x\mc v-\mc u&-&\text{shear displacement,}\\
p_2=\rho\p_t\mc v&-&\text{momentum,}\\
p_3=\p_x \mc u&-&\text{angular displacement,}\\
p_4=I_{\rho}\p_t\mc u&-&\text{angular momentum,}
\end{array}\right.
\end{displaymath}
where we dropped the arguments for clarity of notation, the system \eqref{secondorderbeam} can be written as a first order system for $\mb p=(p_i)_{i=1,\ldots,4}$ of the form
\begin{equation}\label{firstorderbeam}
\partial_t\mb p=\mc M\partial_x\mb p+\mc N\mb p,
\end{equation}
where
\begin{displaymath}
\mc M=\left(\begin{array}{cccc}
0&\frac{1}{\rho}&0&0\\
K&0&0&0\\0&0&0&\frac{1}{I_{\rho}}\\
0&0&EI&0
\end{array}\right),\quad
\mc N=\left(\begin{array}{cccc}
0&0&0&-\frac{1}{I_{\rho}}\\
0&0&0&0\\
0&0&0&0\\
K&0&0&0
\end{array}\right).
\end{displaymath}
The eigenvalues of the matrix $\mc M$ are $\lam_{1,3}=\mp\sqrt{\frac{K}{\rho}}, \lam_{2,4}=\mp\sqrt{\frac{EI}{I_{\rho}}}$, and the corresponding eigenvectors form the diagonalizing matrix
\begin{displaymath}
\mc F=\left(\begin{array}{cccc}
1&0&1&0\\
-\sqrt{K\rho}&0&\sqrt{K\rho}&0\\
0&1&0&1\\
0&-\sqrt{EII_{\rho}}&0&\sqrt{EII_{\rho}}
\end{array}\right).
\end{displaymath}
We can apply our theory to the considered problem by Remark \ref{uwagidoteorii}(ii).  Upon introducing the Riemann invariants $\mb u=\mc F^{-1}\mb p$, the system \eqref{firstorderbeam} is equivalent to
\begin{equation}
\partial_t\mb u=\left(\begin{array}{cccc}
-\sqrt{\frac{K}{\rho}}&0&0&0\\
0&-\sqrt{\frac{EI}{I_{\rho}}}&0&0\\
0&0&\sqrt{\frac{K}{\rho}}&0\\
0&0&0&\sqrt{\frac{EI}{I_{\rho}}}
\end{array}\right)\partial_x\mb u+\mc F^{-1}\mc N\mc F\mb u.
\end{equation}
As noted in Introduction, we will  focus on the principal part of the problem, that is, we assume that $\mc N=0,$ and consider  two connected beams which are represented by the graph depicted in Figure \ref{beamsgraph}. The problem is given by
\begin{equation}\label{twobeams}
\left.\begin{array}{ll}
\partial_t\mb p^1=\mc M^1\partial_x\mb p^1&\text{on $\mb e_1$},\\
\partial_t\mb p^2=\mc M^2\partial_x\mb p^2&\text{on $\mb e_2$}.
\end{array}\right.
\end{equation}
\begin{figure}[ht]
\centering
\begin{tikzpicture}
\draw [-] (-4,0) -- (-0,0);
\draw [-] (0.0,0) -- (4.05,0);

\draw (0,0)node{$\bullet$};
\draw (4,0)node{$\bullet$};
\draw (-4,0)node{$\bullet$};

\draw (-4.0,0.2) node {\footnotesize{$\bv_1$}};
\draw (-0.0,0.2) node {\footnotesize{$\bv_2$}};
\draw (4.0,0.2) node {\footnotesize{$\bv_3$}};
\draw (-2,-0.2)node{$\footnotesize{\mb e_1}$};
\draw (2,-0.2)node{$\footnotesize{\mb e_2}$};
\end{tikzpicture}
\caption{A graph representing two connected beams.}
\label{beamsgraph}
\end{figure}
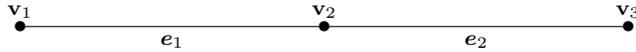
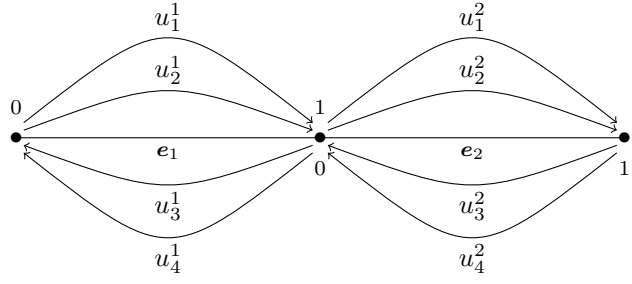
\begin{figure}
\centering
\begin{tikzpicture}
\draw [-] (-4,0) -- (-0,0);
\draw [-] (0.0,0) -- (4.05,0);

\draw (-2,-0.2)node{$\footnotesize{\mb e_1}$};
\draw (2,-0.2)node{$\footnotesize{\mb e_2}$};

\draw (-4.0,0.4) node {\footnotesize{$0$}};
\draw (-0.0,0.4) node {\footnotesize{$1$}};
\draw (-0.0,-0.4) node {\footnotesize{$0$}};
\draw (4.0,-0.4) node {\footnotesize{$1$}};

\draw (0,0)node{$\bullet$};
\draw (4,0)node{$\bullet$};
\draw (-4,0)node{$\bullet$};

\draw (-2,1.6)node{$\footnotesize{u^1_1}$};
\draw (-2,0.9)node{$\footnotesize{u^1_2}$};
\draw (-2,-0.9)node{$\footnotesize{u^1_3}$};
\draw (-2,-1.6)node{$\footnotesize{u^1_4}$};
\draw (2,1.6)node{$\footnotesize{u^2_1}$};
\draw (2,0.9)node{$\footnotesize{u^2_2}$};
\draw (2,-0.9)node{$\footnotesize{u^2_3}$};
\draw (2,-1.6)node{$\footnotesize{u^2_4}$};

\draw[->] (-3.9,0.2)..controls(-2,1.7)..(-0.1,0.2);
\draw[->] (-3.9,0.1)..controls(-2,0.8) ..(-0.1,0.1);
\draw[<-] (-3.9,-0.1)..controls(-2,-0.8) ..(-0.1,-0.1);
\draw[<-] (-3.9,-0.2)..controls(-2,-1.7) ..(-0.1,-0.2);
\draw[<-] (3.9,0.2)..controls(2,1.7)..(0.1,0.2);
\draw[<-] (3.9,0.1)..controls(2,0.8) ..(0.1,0.1);
\draw[->] (3.9,-0.1)..controls(2,-0.8) ..(0.1,-0.1);
\draw[->] (3.9,-0.2)..controls(2,-1.7) ..(0.1,-0.2);

\end{tikzpicture}
\caption{Riemann invariants for the problem \eqref{twobeams}.}
\end{figure}
We parametrize the edges as $l_1(\bv_1)=l_2(\bv_2)=0$ and $l_1(\bv_2)=l_2(\bv_3)=1$. We assume that the beams are made of the same material, that is, $\mc M^1=\mc M^2=\mc M.$  To simplify the calculations, we also assume that the transport velocities equal $1$.  Denoting by $\mc I$ the $2\times 2$ identity matrix and $\msf C={diag}(-\mc I,\mc I,-\mc I,\mc I)$, after diagonalization we obtain the following system:
\begin{equation}\label{beamdiagonalized}
\p_t\binom{\mb u^1}{\mb u^2}=\msf C\p_x\binom{\mb u^1}{\mb u^2}.
\end{equation}
According to the theory of \cite{JBAB1}, the correct numbers of boundary conditions at the vertices are 2 at $\bv_1$ and $\bv_3,$ and 4 at $\bv_2$. We assume that they are as follows.
\begin{equation}
\label{bccond}
\begin{tabular} {|c|c|c|}
\hline
$\bv_1$ & $\bv_2$&$\bv_3$\\
\hline
$\left.\begin{array}{cc}p_1^1(0)=\frac{3}{\sqrt{K\rho}}p_2^1(0)\\p_3^1(0)=\frac{3}{\sqrt{EII_{\rho}}}p_4^1(0) \end{array}\right.$ & $\left.\begin{array}{cccc}p_1^1(1)=p_1^2(0)\\p_2^1(1)=p_2^2(0)\\p_3^1(1)=0\\p_3^2(0)=0 \end{array}\right.$ &$\left.\begin{array}{cc}-3p_1^2(1)=\frac{5}{\sqrt{K\rho}}p_2^2(1)\\-3p_3^2(1)=\frac{5}{\sqrt{EII_{\rho}}}p_4^2(1) \end{array}\right.$ \\
\hline
\end{tabular}
\end{equation}
For the Riemann invariants we have
\begin{displaymath}
\mb\Psi^{out}_{\bv_1}\binom{u_1^1(0)}{u_2^1(0)}=\mb\Psi^{in}_{\bv_1}\binom{u_3^1(0)}{u_4^1(0)},\quad \mb\Psi^{out}_{\bv_1}=\left(\begin{array}{cc} 1&0\\0&1 \end{array}\right),\quad \mb\Psi^{in}_{\bv_1}=\left(\begin{array}{cc} \frac{1}{2}&0\\0&\frac{1}{2} \end{array}\right),
\end{displaymath}
\begin{displaymath}
\mb\Psi^{out}_{\bv_2}\left(\begin{array}{c} u_3^1(1)\\u_4^1(1)\\u_1^2(0)\\u_2^2(0) \end{array}\right)=\mb\Psi^{in}_{\bv_2}\left(\begin{array}{c} u_1^1(1)\\u_2^1(1)\\u_3^2(0)\\u_4^2(0) \end{array}\right),
\end{displaymath}
\begin{displaymath}
\mb\Psi^{out}_{\bv_2}=\left(\begin{array}{cccc}
1&0&-1&0\\
1&0&1&0\\
0&1&0&0\\
0&0&0&1
\end{array}\right),\quad
\mb\Psi^{in}_{\bv_2}=\left(\begin{array}{cccc}
-1&0&1&0\\
1&0&1&0\\
0&-1&0&0\\
0&0&0&-1
\end{array}\right),
\end{displaymath}
\begin{displaymath}
\mb\Psi^{out}_{\bv_3}\binom{u_3^2(1)}{u_4^2(1)}=\mb\Psi^{in}_{\bv_3}\binom{u_1^2(1)}{u_2^2(1)},\quad \mb\Psi^{out}_{\bv_3}=\left(\begin{array}{cc} 1&0\\0&1 \end{array}\right),\quad \mb\Psi^{in}_{\bv_3}=\left(\begin{array}{cc} \frac{1}{4}&0\\0&\frac{1}{4} \end{array}\right).
\end{displaymath}
To write our problem as a port-Hamiltonian system \eqref{syst1}, let us introduce
$$
\mb\upsilon=(u_1^1,u_2^1,u_1^2,u_2^2),\quad \mb\w=(u_3^1,u_4^1,u_3^2,u_4^2)
$$
and let $\mc P$ be the permutation matrix interchanging rows $5$ with $3$ and $6$ with $4$.
Then, since $(\mb\ups,\mb\w)^T=\mc P(\mb u^1,\mb u^2)^T$, \eqref{beamdiagonalized} becomes
$$
\p_t\binom{\mb\ups}{\mb\w}=
\left(\begin{array}{cc}-\mc C_+&0\\0&\mc C_-  \end{array}\right)\p_x\binom{\mb\ups}{\mb\w},
$$
where $\mc C_+=\mc C_-={diag}(\mc I,\mc I)$ and the matrices in \eqref{syst1bc} are given by
\begin{displaymath}
\mb\Xi_{out}=\left(\begin{array}{cccccccc}
1&0&0&0&0&0&0&0\\
0&1&0&0&0&0&0&0\\
0&0&-1&0&1&0&0&0\\
0&0&1&0&1&0&0&0\\
0&0&0&0&0&1&0&0\\
0&0&0&1&0&0&0&0\\
0&0&0&0&0&0&1&0\\
0&0&0&0&0&0&0&1\\
\end{array}\right),
\quad
\mb\Xi_{in}=\left(\begin{array}{cccccccc}
0&0&0&0&\frac{1}{2}&0&0&0\\
0&0&0&0&0&\frac{1}{2}&0&0\\
-1&0&0&0&0&0&1&0\\
1&0&0&0&0&0&1&0\\
0&-1&0&0&0&0&0&0\\
0&0&0&0&0&0&0&-1\\
0&0&\frac{1}{4}&0&0&0&0&0\\
0&0&0&\frac{1}{4}&0&0&0&0\\
\end{array}\right),
\end{displaymath}
so that the boundary matrix is given by
\begin{displaymath}
\mc B=\mb\Xi_{out}^{-1}\mb\Xi_{in}=\left(\begin{array}{cccccccc}
0&0&0&0&\frac{1}{2}&0&0&0\\
0&0&0&0&0&\frac{1}{2}&0&0\\
1&0&0&0&0&0&0&0\\
0&0&0&0&0&0&0&-1\\
0&0&0&0&0&0&1&0\\
0&-1&0&0&0&0&0&0\\
0&0&\frac{1}{4}&0&0&0&0&0\\
0&0&0&\frac{1}{4}&0&0&0&0\\
\end{array}\right).
\end{displaymath}
Routine calculations show that the characteristic polynomial is given by
$$
\left(\la^4 -\frac{1}{8}\right)\left(\la^2+\frac{1}{2}\right) \left(\la^2+\frac{1}{4}\right)= 0
$$
so that $\max\{|\la|\;:\;\la \in \sigma(\mc B)\}= \frac{1}{\sqrt{2}}$ and right and left eigenvectors are given by
\begin{center}
\begin{tabular}{|c|c|c|c|}
\hline
eigenvalues &$\la^4 -\frac{1}{8}=0$&$\left(\beta^2+\frac{1}{2}\right)=0$&$\left(\gamma^2+\frac{1}{4}\right)= 0$\\
\hline
right eigenvectors&$(1,0,\frac{1}{\la},0,2\la,0,2\la^2,0)$&$(0,1,0,0,0,2\beta,0,0)$&$(0,0,0,1,0,0,0,-\gamma)$\\
\hline
left eigenvectors&$(1,0,\la,0,\frac{1}{2\la},0,4\la^2,0)$&$(0,1,0,0,0,-\beta,0,0)$&$(0,0,0,1,0,0,0,4\gamma)$\\
\hline\end{tabular}\,.
\end{center}
 If we denote by $\mb F_i$ and $\mb E_i$ the left and right eigenvectors belonging to $i=\la,\beta, \gamma,$
 we see that
 $$
 \mb F_\la\cdot \mb E_\la =4, \quad \mb F_\beta\cdot \mb E_\beta =2, \quad\mb F_\gamma\cdot \mb E_\gamma =2.
 $$
 Hence, by \eqref{Urep}, we can estimate the norm of projectors in $\mbb C^8$ with $|\cdot|_1$ norm as follows:
 $$
 \|\Pi_\la\mb u\|_1 \leq \frac{|\mb F_\la\mb  u|}{\mb F_\la\cdot\mb E_\la}\|\mb E_\la\|_1 \leq \frac{1}{4}\|\mb F_\la\|_\infty\|\mb E_\la\|_1\|\mb u\|_1\approx 1.62 \|\mb u\|_1.
 $$
 Similarly,
 \begin{align*}
 \|\Pi_\beta\mb u\|_1 & \leq \frac{1}{2}\|\mb F_\beta\|_\infty\|\mb E_\beta\|_1\|\mb u\|_1\approx 1.21 \|\mb u\|_1,\\
 \|\Pi_\gamma\mb u\|_1 & \leq \frac{1}{2}\|\mb F_\gamma\|_\infty\|\mb E_\gamma\|_1\|\mb u\|_1=1.5 \|\mb u\|_1.
 \end{align*}
 Denote by $\sem{G_0}$ the semigroup generated by \eqref{twobeams} with boundary conditions \eqref{bccond}. Then, as in the calculations in the proof of Lemma \ref{aslemma}, in the operator norm of $\mb X_1$ we have
 $$
\|G_0(t)\|\leq M e^{-\frac{\ln 2}{2} t}, \quad t\geq 0,
$$
 where $M = \la (4 \cdot 1.62 + 2\cdot 1.21 + 2\cdot 1.5)\approx 8.41$ and $\lambda$ can be taken $\max\{|\lambda|:\lambda\in\sigma(\mc B)\}$ since all eigenvalues are simple.
  \end{example}

\def\cprime{$'$}

 \end{document}